\documentclass{article}
\def\disc{{\rm disc}}
\def\SL{{\bf{\rm SL}}}
\def\Q{\hbox{\bf Q}}
\def\Tr{\mathop{\rm Tr}\nolimits}
\def\Z{\hbox{\bf Z}}
\def\mod{\mathop{\;\rm mod}\nolimits}

\def\F{\hbox{\bf F}}

\begin{document}

\begin{center}
{\bf Extensions g\'en\'eriques de groupe de Galois $\SL_2(\F_3)$}
\\ Jean-Fran\c cois Mestre.
\end{center}

\medskip 

\section{Introduction}
Soit $G$ un groupe fini. Un polyn\^ome g\'en\'erique (sur $\Q$) pour $G$ est un 
\'el\'ement de $\Q(t_1,\ldots,t_s)[X]$, o\`u les $t_i$ sont des
ind\'etermin\'ees, de groupe de Galois $G$ sur $\Q(t_1,\ldots,t_s)$, et tel
que, si $K$ est un corps contenant $\Q$ et $L/K$ une extension de groupe de
Galois contenue dans $G$, il existe une sp\'ecialisation des param\`etres
$t_1,\ldots,t_s$ en des valeurs de $K$ telle que le corps des racines
du polyn\^ome $P$ ainsi obtenu soit \'egal \`a $L$.

Le groupe de plus petit ordre pour lequel on ignore s'il existe ou non un 
polyn\^ome
g\'en\'erique est le groupe $\SL_2(\F_3)=L_3(2).$ 

\medskip
Nous prouvons ici:

\medskip

{\sc Th\'eor\`eme 1.-} {\it
Il existe un polyn\^ome g\'en\'erique pour le groupe de
Galois $\SL_2(\F_3)$.}

\medskip

Pour ce faire, nous prouvons la proposition suivante :

{\sc Proposition 1.-} {\it Le polyn\^ome 
$$P=X^4-4X^3+36VX^2+2U^2X^2+4U^2X-8U^2VX+36U^2V^2+U^4-4U^2V$$
est un polyn\^ome g\'en\'erique pour les extensions
de groupe de Galois $A_4$.

De plus, l'invariant de Witt de la forme quadratique $\Tr(x^2)$
associ\'ee est \'egal \`a
$$(-1,-1)+(U^2-9,-2(U^2V-9V+1-U^2)).$$}

Le th\'eor\`eme en est une cons\'equence facile : le groupe
$\SL_2(\F_3)$ est isomorphe au groupe $\tilde{A}_4$. Or,
si $K/k$ est une extension de degr\'e $4$ dont la cl\^oture galoisienne $L$ 
a comme
groupe de Galois $A_4$, l'obstruction
\`a l'existence d'une extension quadratique $M/L$ telle que $M/k$
soit galoisienne de groupe de Galois $\tilde{A}_4$ est \'egale
\`a l'invariant de Witt de la forme $\Tr_{K/k}(x^2)$. 

Or, si $(a,b)\in k^*$, $k$ corps de caract\'eristique diff\'erente de $2$,
il est bien connu que 
le symbole $(-1,-1)+(a,b)$ est nul sur $k$ si et seulement s'il   
existe $A,B,C,D,E\in k$ tels que

$$a=D^2(1+A^2+A^2B^2)(1+B^2+B^2C^2),\;\;\;b=E^2(1+B^2+B^2C^2)(1+C^2+C^2A^2).$$

En posant $a=(U-3)/(U+3)$ et $b=-2(U^2V-9V+1-U^2)$, on en d\'eduit donc
imm\'ediatement une param\'etrisation de $(U,V)$ en fonction des $5$
param\`etres $A,B,C,D,E$.

\medskip

{\sc Remarque.-} Pour obtenir un
un polyn\^ome g\'en\'erique explicite 
pour $L_3(2)$, on peut par exemple appliquer les r\'esultats de
T. Crespo (Explicit construction of $2 S_n$ Galois extensions, J. Algebra 129 (1990)).

\section{D\'emonstration de la proposition $1$.}

Pour tout polyn\^ome $f$ en une ind\'etermin\'ee $X$, on note $\disc(f)$
son discriminant.

Dans ce qui suit, on note $k$ le corps $\Q(a_1,a_2,a_3,a_4)$, 
o\`u les $a_i$ sont des ind\'etermin\'ees. 

On d\'emontre ais\'ement la proposition suivante:

\medskip

{\sc Proposition $2$.-} {\it Soit $P$ l'\'el\'ement de $k[X]$ \'egal \`a
$X^4+a_1X^3+a_2X^2+a_3X+a_4$, et
$Q=b_0X^3+b_1X^2+b_2X+b_3$ le reste de la division euclidienne de
$P'^2$ par $P$. On a :

1) Si $x_1,\ldots,x_4$ sont les racines de $P$, les racines de $Q$ sont $$\frac{x_1x_2-x_3x_4}{x_1+x_2-x_3-x_4}$$ et les
deux autres quantit\'es analogues. (Donc $Q$ est une r\'esolvante pour le groupe
$(\Z/2\Z)^2.$)

2) Soit $T$ une ind\'etermin\'ee, et $P_T=P-TQ\in k(T)[X].$ 
Il existe un polyn\^ome $U(T)\in k[T]$ tel que le reste de $(P'_T)^2\mod P$ est \'egal \`a $U(T) Q$.

3)  On a
$\disc(P_T)=\disc(P) U(T)^2$. 

4) L'extension $M$ de $k(T)$ obtenue par adjonction des racines de $P_T$ 
contient
le corps $L$ des racines de $Q$. L'extension $M/L(T)$ est une extension
r\'eguli\`ere de groupe de Galois $(\Z/2\Z)^2$.

5) \`A constante multiplicative pr\`es, le polyn\^ome $Q$ est l'unique
polyn\^ome du troisi\`eme degr\'e tel que le groupe de Galois de
$P-TQ$ sur $\overline{k}(T)$ soit \'egal \`a $(\Z/2\Z)^2.$

6) Soit $\tilde{P}=Y^4P(X/Y)$ le polyn\^ome homog\'en\'eis\'e de $P$,
et $H=P"_{X^2}P"_{Y^2}-P"_{XY}^2$ le Hessien de $P$. On a 
$H(X,1)=(-9a_1^2+24a_2)P-9Q.$

7) On a $\disc(Q)=\disc(P) b_0^2$.
De plus, l'application qui \`a $(a_1,a_2,a_3,a_4)$ associe $(a_1,c_1,c_2,c_3)$,
o\`u $c_i=b_i/b_0$, $1\leq i\leq 3$, est birationnelle.}
\medskip

En fait, seul le point $7)$ nous sera utile dans la suite. 
Explicitons les formules donnant les $b_i$ en fonction des $a_i$:

On a $b_3=a_3^2-a_1^2a_4,b_2=4a_2a_3-a_1^2a_3-8a_1
a_4,b_1=-a_1^2a_2-2a_1a_3-16a_4+4a_2^2,b_0=-8a_3-a_1^3+4a_1a_2,$
d'o\`u
$$a_2=c_1a_1-2c_2,a_3=c_2a_1-8c_3, a_4 = c_3a_1+c_2^2-4c_1c_3.$$

\medskip
La param\'etrisation des polyn\^omes de degr\'e $4$ de  discriminant carr\'e
contenu dans $A_4$ est ainsi ramen\'ee \`a celle des polyn\^omes
de degr\'e $3$ de discriminant carr\'e. 

On peut par exemple prendre

$$c_3 = 1/27c_1^3-1/12c_1v^2-27/4u^3-1/4uv^2-9/4c_1u^2, c_2 = c_1^2/3-27u^2/4-v^2/4,$$
d'o\`u $a_1,a_2,a_3,a_4$ en fonction de $a_1,c_1,u,v$. 

On trouve alors le polyn\^ome de la proposition $1$ par le changement de
va\-ria\-bles 
$$\left\{
\begin{array}{l}
c_1=\frac{3a_1V-18uV+6u}{4V}\\
v=\frac{uU}{V}\\
x=\frac{-2u-Va_1+2Xu}{4V}
\end{array}\right.$$
qui transforme birationnellement $(c_1,v,x)$ en $(U,V,X)$,
l'application r\'eciproque \'etant
$$\left\{
\begin{array}{l}
U = -\frac{6v}{-18u-4c_1+3a_1} \\
V = -\frac{6u}{-18u-4c_1+3a_1}\\
X = -\frac{2(6x+9u+2c_1)}{-18u-4c_1+3a_1}
\end{array}\right.$$

Un calcul sans
difficult\'e permet alors de montrer que l'invariant de Witt de la forme
$\Tr(x^2)$ correspondante a la forme voulue, d'o\`u le r\'esultat.  

\end{document}